\documentclass[11]{article}
\usepackage[english]{babel}
\usepackage{amssymb}
\usepackage{inputenc}
\usepackage{graphicx}
\begin{document}
\large \noindent
\begin{center}
\textbf{\Large Optimal quadrature formulas of closed type
in the space  $L_2^{(m)}(0,1)$}
\\
\textbf{Kh.M.Shadimetov}\footnote{This article is the fifth
section of the author's Candidate dissertation entitled "Optimal
formulas of approximate integration for
differentiable functions"\ - Novosibirsk, 1983. 140 p.}\\
{\it Institute of mathematics and information technologies,\\
Uzbek Academy of Sciences, Tashkent, Uzbekistan}
\end{center}

\begin{abstract}
It is discussed the problem on construction of optimal quadrature
formulas in the sense of Sard in the space $L_2^{(m)}(0,1)$, when
the nodes of quadrature formulas are equally spaced. Here the
representations of optimal coefficients for any natural numbers
$m$ and $N$ are found.
\end{abstract}

\textbf{2000 Mathematics Subject Classification:} 65D32\\

\textbf{Key words and phrases:} \emph{Sobolev space, optimal
quadrature formula, positive coefficients, error functional.}\\

Consider quadrature formulas of the form
    $$
\int\limits_0^1 {\varphi (x)dx \cong \sum\limits_{\beta  = 0}^N
{C[\beta ]\,\varphi [\beta ]} }\eqno (1)
$$
with the \emph{error functional}
$$
\ell(x)=\varepsilon_{[0,1]}(x)-\sum\limits_{\beta=0}^NC[\beta]\delta(x-h\beta),\eqno
(2)
$$
where $C[\beta]$ are the coefficients of formula (1),
$[\beta]=h\beta$, $h=1/N$, $N=1,2,3,...$, the function
$\varphi(x)$ belongs to the space $L_2^{(m)}(0,1)$. Norm of
functions in the space $L_2^{(m)}(0,1)$ is determined by formula
$$
\|\varphi(x)|L_2^{(m)}(0,1)\|=\left(\int\limits_0^1
\left(\varphi^{(m)}(x)\right)^2dx\right)^{\frac{1}{2}}.
$$
Below $[\beta]=h\beta$.

Quadrature formulas of the form  (1) is \emph{closed type}, since
the points $x=0$ and $x=1$ are the nodes of the formula.

The difference
    $$
(\ell(x),\varphi(x))=\int\limits_0^1 {\varphi (x)dx -
\sum\limits_{\beta = 0}^N {C[\beta ]\,\varphi [\beta ]} }
$$
is called \emph{the error} of formula (1).

\emph{The problem} of construction of optimal quadrature formulas
(1) in the sense of Sard in the space $L_2^{(m)}(0,1)$ consists of
computation   the quantity
$$
\|\stackrel{\circ}{\ell}(x)|L_2^{(m)*}(0,1)\|=\inf_{C[\beta]}\sup_{\|\varphi(x)\|\neq
0}\frac{|(\ell(x),\varphi(x))|}{\|\varphi(x)|L_2^{(m)}(0,1)\|}.\eqno
(3)
$$
where $L_2^{(m)*}(0,1)$ is the conjugate space to the space
$L_2^{(m)}(0,1)$.

The coefficients $C[\beta]$ satisfying equality (3) (if there
exist) are called \emph{optimal} and are denoted as
$\stackrel{\circ}{C}[\beta]$. Formulas of the form (1) with
coefficients $\stackrel{\circ}{C}[\beta]$ are called \emph{optimal
quadrature formulas in the sense of Sard}.

The problem of construction of optimal quadrature formulas in the
sense of Sard in different spaces are investigated by many
mathematicians. In the space $L_2^{(m)}$ this problem was
investigated by S.L.Sobolev in [1] and  for optimal coefficients
$\stackrel{\circ}{C}[\beta]$ the system of linear equations was
obtained, which for $x_{\beta}=[\beta]=h\beta$
$(\beta=\overline{0,N})$ has the form
$$
\sum\limits_{\gamma  = 0}^N {\stackrel{\circ}{C}[\gamma]
{{|[\beta] - [\gamma] |^{2m - 1} } \over {2 \cdot (2m - 1)!}} +
\sum\limits_{\alpha  = 0}^{m - 1} {\lambda _\alpha  [\beta]
^\alpha   = {{[\beta]^{2m} + (1 - [\beta])^{2m} } \over {2\left(
{2m} \right)!}},\ \ \ \ [\beta] \in [0,1],} } \eqno (4)
$$
$$
\sum\limits_{\gamma  = 0}^N {\stackrel{\circ}{C}[\gamma] [\gamma]
^\alpha   = {1 \over {\alpha  + 1}},\,\,\,\,\,\alpha  = \overline
{0,m - 1}}. \eqno (5)
$$
In this system $\stackrel{\circ}{C}[\beta]$ and
$\lambda_{\alpha}$, $\alpha=0,1,...,m-1$ are unknowns.

Note that the existence and uniqueness of the solution of the
system (4)-(5) was proved by S.L.Sobolev in [1].

The aim of this work is finding the explicit forms of the  optimal
coefficients $\stackrel{\circ}{C}[\beta]$.

Later for convenience optimal coefficients
$\stackrel{\circ}{C}[\beta]$ we remain as $C[\beta]$.

Here mainly is used the concept of discrete argument functions
and operations on them. For completeness we give some definitions
about functions of discrete argument.

Assume that $\varphi(x)$ and $\psi(x)$ are real-valued functions
of real variable and are defined in real line $\mathbb{R}$.

\textbf{Definition 1.} Function $\varphi(h\beta )$ is called
\emph{the function of discrete argument}, if it is given on some
set of integer values of $\beta $.

\textbf{Definition 2.} \emph{By inner product} of two discrete
functions $\varphi (h\beta )$ and $\psi (h\beta )$ is called the
number
    $$
\left[ {\varphi ,\psi } \right] = \sum\limits_{\beta  =  - \infty
}^\infty  {\varphi (h\beta ) \cdot \psi (h\beta )},
$$
if the series on right hand side of the last equality converges
absolutely.

\textbf{Definition 3.} \textit{By convolution} of two functions
$\varphi (h\beta )$ and $\psi (h\beta )$ is called the inner
product
$$\varphi (h\beta )*\psi (h\beta ) = \left[ {\varphi (h\gamma
),\psi (h\beta  - h\gamma )} \right] = \sum\limits_{\gamma  =  -
\infty }^\infty  {\varphi (h\gamma ) \cdot \psi (h\beta  - h\gamma
)}.
$$

Suppose $C[\beta]=0$ for $\beta<0$ and $\beta>N$. Then, using
definition 3, system (4)-(5) can be rewritten in the following
form
$$
G_{m,1}[\beta]*C[\beta]+P_{m-1}[\beta]=f_m[\beta]\mbox{ when }
[\beta]\in [0,1],\eqno (6)
$$
$$
C[\beta]=0 \mbox{ when } [\beta]\not\in [0,1], \eqno (7)
$$
$$
\sum\limits_{\beta=0}^{N}C[\beta]\
[\beta]^{\alpha}=\frac{1}{\alpha+1},\ \ \alpha=0,1,2,...,m-1,
\eqno (8)
$$
where $P_{m-1}[\beta]$ is polynomial of degree $m-1$ with respect
to $[\beta]$ and
$$
G_{m,1}[\beta]=\frac{[\beta]^{2m-1}\mathrm{sign}[\beta]}{2
(2m-1)!},\eqno (9)
$$
$$
f_m[\beta]=\frac{[\beta]^{2m}+(1-[\beta])^{2m}}{2(2m)!}.\eqno (10)
$$

Consider the following problem.

\textbf{Problem A.} \emph{Find the discrete function $C[\beta]$
and unknown polynomial $P_{m-1}[\beta]$.}

Denote
$$
v[\beta]=G_{m,1}[\beta]*C[\beta] \eqno (11)
$$
and
$$
u[\beta]=v[\beta]+P_{m-1}[\beta]. \eqno (12)
$$
Then we have to express $C[\beta]$ with the help of $u[\beta]$.
For this we must construct discrete operator $D_m[\beta]$, which
satisfies the equation
$$
hD_m[\beta]*G_{m,1}[\beta]=\delta[\beta], \eqno (13)
$$
where $\delta[\beta]$ equals 0 when $\beta\neq 0$ and equals 1
when $\beta=0$, i.e. $\delta[\beta]$ is the discrete
delta-function, $G_{m,1}[\beta]$ is defined by formula (9).

In connection with this in the work [2] the discrete analogue
$D_m[\beta]$ of the differential operator $d^{2m}/dx^{2m}$ was
constructed.

In [2] the following was proved.

{\bf Theorem 1.} {\it The discrete analogue of the differential
operator ${d^{2m}}/{dx^{2m}}$ have the form
$$
D_m[\beta]=\frac{(2m-1)!}{h^{2m}}\left\{
\begin{array}{lll}
{\displaystyle
\sum\limits_{k=1}^{m-1}\frac{(1-q_k)^{2m+1}q_k^{|\beta|}}
{q_kE_{2m-1}(q_k)}        }
& \mbox{ for }& |\beta|\geq 2,\\
{\displaystyle 1+\sum\limits_{k=1}^{m-1}\frac{(1-q_k)^{2m+1}}
{E_{2m-1}(q_k)}                  }
& \mbox{ for }& |\beta|= 1,\\
{\displaystyle
-2^{2m-1}+\sum\limits_{k=1}^{m-1}\frac{(1-q_k)^{2m+1}} {q_k
E_{2m-1}(q_k)}           } & \mbox{ for }& \beta= 0,
\end{array}
\right. \eqno (14)
$$
where $E_{2m-1}(x)$ is the Euler-Frobenius polynomial of degree
$2m-1$, $q_k$ are roots of the polynomial $E_{2m-2}(x)$,
$|q_k|<1$, $h$ is small parameter.}

{\bf Theorem 2.} {\it For the operator $D_m[\beta]$ and monomials
$[\beta]^k=(h\beta)^k$ the following are true}
$$
\sum_{\beta}D_m[\beta][\beta]^k= \left\{
\begin{array}{lll}
0&\mbox{ \emph{when} } & 0\leq k\leq 2m-1,\\
(2m)!&\mbox{ \emph{when} } & k= 2m,
\end{array}
\right. \eqno (15)
$$
$$
\sum_{\beta}D_m[\beta][\beta]^k= \left\{
\begin{array}{lll}
0&\mbox{ \emph{when} } & 2m+1\leq k\leq 4m-1,\\
{\displaystyle             \frac{h^{2m}(4m)!B_{2m}}{(2m)!}    }
&\mbox{ \emph{when} }  & k= 4m.
\end{array}
\right. \eqno (16)
$$

Taking into account (13) and theorems 1, 2 for optimal
coefficients we have
$$
C[\beta]=hD_m[\beta]*u[\beta]. \eqno (17)
$$

Thus, if we find the function $u[\beta]$, then optimal
coefficients will be found from (17).

In order to calculate convolution (17) we need representation of
the function $u[\beta]$ in all integer values of $\beta$. From (6)
we have $u[\beta]=f_m[\beta]$ when $[\beta]\in [0,1]$. Now we must
find representation of $u[\beta]$ for $\beta<0$ and $\beta>N$.

Since $C[\beta]=0$ when $[\beta]\not\in [0,1]$, then
$$
C[\beta]=hD_m[\beta]*u[\beta]=0,\ \ \ [\beta]\not\in [0,1].
$$

We calculate the convolution $v[\beta]=G_{m,1}[\beta]*C[\beta]$
when $[\beta]\not\in [0,1]$.

Suppose $\beta<0$, then taking into account (4), (5) we have
$$
v[\beta]=G_{m,1}[\beta]*C[\beta]=-\sum\limits_{\gamma=-\infty}^{\infty}
C[\gamma]\frac{([\beta]-[\gamma])^{2m-1}}{2(2m-1)!}=
$$
$$
-\frac{1}{2}\sum\limits_{j=0}^{m-1}\frac{[\beta]^{2m-1-j}(-1)^j}{(j+1)!(2m-1-j)!}-
\frac{1}{2}\sum\limits_{j=m}^{2m-1}\frac{[\beta]^{2m-1-j}(-1)^j}{j!(2m-1-j)!}
\sum\limits_{\gamma=0}^NC[\gamma][\gamma]^j.\eqno (18)
$$
Denote
$$
Q_{2m-1}[\beta]=\frac{1}{2}\sum\limits_{j=0}^{m-1}
\frac{[\beta]^{2m-1-j}(-1)^j}{(j+1)!(2m-1-j)!}, \ \
R_{m-1}[\beta]=
\frac{1}{2}\sum\limits_{j=m}^{2m-1}\frac{[\beta]^{2m-1-j}(-1)^j}{j!(2m-1-j)!}
\sum\limits_{\gamma=0}^NC[\gamma][\gamma]^j.
$$
Then from (18) for $v[\beta]$ when $\beta<0$ we obtain
$$
v[\beta]=-Q_{2m-1}[\beta]-R_{m-1}[\beta].\eqno (19)
$$
Similarly for the case $\beta>N$ we have
$$
v[\beta]=Q_{2m-1}[\beta]+R_{m-1}[\beta].\eqno (20)
$$
Denoting
$$
R_{m-1}^-[\beta]=P_{m-1}[\beta]-R_{m-1}[\beta],\eqno (21)
$$
$$
R_{m-1}^+[\beta]=P_{m-1}[\beta]+R_{m-1}[\beta],\eqno (22)
$$
and taking into account (19), (20), (12) we get following problem.

\textbf{Problem B.} {\it Find solution of the equation
$$
hD_m[\beta]*u[\beta]=0,\ \ \  [\beta]\not\in [0,1]\eqno (23)
$$
having the form
$$
u[\beta]= \left\{
\begin{array}{ll}
-Q_{2m-1}[\beta]+R_{m-1}^-[\beta], &\beta<0,\\
f_m[\beta], &0\leq \beta\leq N-1,\\
Q_{2m-1}[\beta]+R_{m-1}^+[\beta], &\beta>N,\\
\end{array}
\right. \eqno (24)
$$
where $R_{m-1}^-[\beta]$ and $R_{m-1}^+[\beta]$ unknown
polynomials of degree $m-1$.}

If we find $R_{m-1}^-[\beta]$ and $R_{m-1}^+[\beta]$, then from
(21), (22) we get
$$
P_{m-1}[\beta]=\frac{1}{2}\left(R_{m-1}^+[\beta]+R_{m-1}^-[\beta]
\right),
$$
$$
R_{m-1}[\beta]=\frac{1}{2}\left(R_{m-1}^+[\beta]-R_{m-1}^-[\beta]
\right).
$$
Unknowns $R_{m-1}^-[\beta]$ and $R_{m-1}^+[\beta]$ can be found
from (23), using the function $D_m[\beta]$. Then we obtain
explicit form of $u[\beta]$ and from (17) will be found the
optimal coefficients $C[\beta]$. Thus, the problem B and
respectively the problem A will be solved.

But here we will not find $R_{m-1}^-[\beta]$ and
$R_{m-1}^+[\beta]$, instead of them, using $D_m[\beta]$ and taking
into account (17), we will find expressions for the optimal
coefficients $C[\beta]$.

The main result of the present work is the following.

\textbf{Theorem 3.} {\it Optimal coefficients of quadrature
formulas of the form (1) on the space $L_2^{(m)}(R)$ have the form
$$
C[\beta]=h\left\{
\begin{array}{ll}
\frac{1}{2}-\sum\limits_{k=1}^{m-1}d_k\ \frac{q_k-q_k^N}{1-q_k}&
\mbox{ for }\beta=0,N,\\
1+\sum\limits_{k=1}^{m-1}
d_k\left(q_k^{\beta}+q_k^{N-\beta}\right)&\mbox{
for } \beta=1,2,...,N-1,\\
0&\mbox{ for }\beta<0, \beta>N,
\end{array}
\right.
$$
where}
$$
\sum\limits_{k=1}^{m-1}d_k\sum_{i=1}^j\frac{q_k+(-1)^{i+1}q_k^{N+i}}
{(q_k-1)^{i+1}}\Delta^i0^j=\frac{B_{j+1}}{j+1},\ \ j=1,2,...,m-1.
$$
\emph{$q_k$ are roots of the Euler-Frobenius polynomial
$E_{2m-2}(x)$ of degree $2m-2$.}

In the proof of theorem  3 the following preliminary results are
used.

The following formula is true [3]
$$
\sum_{\gamma=0}^{n-1}q^{\gamma}\gamma^k=\frac{1}{1-q}\sum\limits_{i=0}^k\left(\frac{q}{1-q}\right)^i
\Delta^i0^k-\frac{q^n}{1-q}\sum_{\gamma=0}^k\left(\frac{q}{1-q}\right)^i
\Delta^i\gamma^k|_{\gamma=n},\eqno (25)
$$
where $\Delta^i\gamma^n$ is finite difference of order $i$ from
$\gamma^n$, $\Delta^i0^k=\Delta^i\gamma^n|_{\gamma=0}$.

And also we use the following well-known formula  (see, for
example, [4])
$$
\sum\limits_{\gamma=0}^{\beta-1}\gamma^k=\sum\limits_{j=1}^{k+1}
\frac{k!B_{k+1-j}}{j!(k+1-j)!}\beta^j.\eqno (26)
$$

It is known in [5], that the Euler-Frobenius polynomials $E_k(x)$
have the form
$$
x E_k(x)=(1-x)^{k+2}D^k\frac{x}{(1-x)^2}, \eqno (27)
$$
where
$$
D=x\frac{d}{dx},\ \ D^k=x\frac{d}{dx}D^{k-1}.
$$
In [5] it was shown that all roots  $q_j^{(k)}$ of the
Euler-Frobenius polynomials $E_k(x)$ are real, negative and
distinct:
$$
q_1^{(k)}<q_2^{(k)}<...<q_k^{(k)}<0. \eqno (28)
$$
Furthermore for the roots (28) the following is true:
$$
q_j^{(k)}\cdot q_{k+1-j}^{(k)}=1.
$$
If we denote $E_k(x)=\sum\limits_{s=0}^ka_s^{(k)}x^s$, then
coefficients $a_s^{(k)}$ of the Euler-Frobenius polynomial
$E_k(x)$ are expressed by formula
$$
a_s^{(k)}=\sum\limits_{j=0}^s(-1)^j{k+2\choose j} (s+1-j)^{k+1}.
$$
This formula was obtained by Euler.

From definition of $E_k(x)$ we get following statements.

{\bf Lemma 1.} {\it For the polynomials $E_k(x)$ the following
recurrence relation holds
$$
E_k(x)=(kx+1)E_{k-1}(x)+x(1-x)E_{k-1}'(x), \eqno (29)
$$
where $E_0(x)=1$, $k=1,2,....$ }

{\bf Lemma 2.} {\it The polynomial $E_k(x)$ satisfies the identity
$$
E_k(x)=x^kE_k\left(\frac{1}{x}\right) \eqno (30)
$$
or otherwise $a_s^{(k)}=a_{k-s}^{(k)},\ \ s=0,1,2,...,k$. }

\textbf{Lemma 3.} {\it Polynomials
$$
P_k(x)=(x-1)^{k+1}\sum\limits_{i=1}^{k+1}\frac{\Delta^i0^{k+1}}{(x-1)^i}\eqno
(31)
$$
and
$$
P_k\left(\frac{1}{x}\right)=\left(\frac{1}{x}-1\right)^{k+1}\sum\limits_{i=1}^{k+1}
\left(\frac{x}{1-x}\right)^i\Delta^i0^{k+1}\eqno (32)
$$
are the Euler-Frobenius polynomials $E_k(x)$ and
$E_k\left(\frac{1}{x}\right)$ respectively. }

\textbf{Lemma 4.} {\it The operator $D_m[\beta]$ satisfies the
relation}
$$
[D_m[\beta],[\beta]^{2m}]=
\sum\limits_{\beta=-\infty}^{\infty}D_m[\beta][\beta]^{2m}=(2m)!.
\eqno (33)
$$

\textbf{Proofs of Lemmas.}

\textbf{Proof of lemma 1.} From (27) obviously that
$$
E_{k-1}(x)=x^{-1}(1-x)^{k+1}D^{k-1}\frac{x}{(1-x)^2}. \eqno (34)
$$
Differentiating by $x$ the polynomial $E_{k-1}(x)$, we obtain
$$
E_{k-1}'(x)=-(1-x)^kx^{-2}(kx+1)D^{k-1}\frac{x}
{(1-x)^2}+\frac{E_k(x)}{x(1-x)}.
$$
Hence and from (34) implies that
$$
(kx+1)E_{k-1}(x)+x(1-x)E_{k-1}'(x)=(kx+1)
x^{-1}(1-x)^{k+1}D^{k-1}\frac{x}{(1-x)^2}-
$$
$$
-(1-x)^{k+1}x^{-1}(kx+1)D^{k-1}\frac{x}{(1-x)^2}+ E_k(x)=E_k(x).
$$
Thus lemma 1 is proved.

\textbf{Proof of lemma 2.} Lemma 2 we will proof by induction. For
$k=1$ from (27) we find
$$
E_1(x)=x+1.
$$
Assume that for $k\geq 1$ the equality
$a_n^{(k-1)}=a_{k-1-n}^{(k-1)}$, $n=0,1,...,k-1$ is fulfilled.
Suppose that $a_n^{(k-1)}=0$ for  $n<0$ and $n>k-1$.

From (29) we get
 $$
 a_s^{(k)}=(s+1)a_s^{(k-1)}+(k-s+1)a_{s-1}^{(k-1)};
 $$
then, using assumptions of the induction, we obtain
 $$
 a_{k-s}^{(k)}=(k-s+1)a_{k-s}^{(k-1)}+(s+1)a_{k-s-1}^{(k-1)}=
 (k-s+1)a_{s-1}^{(k-1)}+(s+1)a_s^{(k-1)}=a_s^{(k)},
 $$
and Lemma 2 is proved.

\textbf{Proof of lemma 3.} Consider relations (29)
$$
E_k(x)=(kx+1)E_{k-1}(x)+x(1-x)E_{k-1}'(x),
$$
$$
E_0(x)=1.
$$
If $P_k(x)$ also satisfies this relations then lemma 3 will be
proved. From (31) evidently that $P_0(x)=1$. We denote
$$
V_k(x)=(kx+1)P_{k-1}(x)+x(1-x)P_{k-1}'(x).
$$
Using the equalities
$$
P_{k-1}(x)=\sum_{i=1}^k\frac{(x-1)^k}{(x-1)^i}\Delta^i0^k,
$$
$$
P_{k-1}'(x)=k(x-1)^{k-1}\sum_{i=1}^k\frac{\Delta^i0^k}{(x-1)^i}-
(x-1)^k\sum_{i=1}^k\frac{i\Delta^i0^k}{(x-1)^{i+1}},
$$
after some simplifications we have
$$
V_k(x)=(kx+1)P_{k-1}(x)+x(1-x)P_{k-1}'(x)=
$$
$$
=P_{k-1}(x)+x(x-1)^k\sum_{i=1}^k\frac{\Delta^i0^k}{(x-1)^i}=
\sum_{i=1}^k\Delta^i0^k(1+xi)(x-1)^{k-i}.
$$
Doing change of variables $x-1=y$, we obtain
$$
V_k(y+1)=\sum_{i=0}^ky^{k-i}(\Delta^i0^k+\Delta^{i+1}0^k)(i+1),
$$
$$
P_k(y+1)=\sum_{i=0}^ky^{k-i}\Delta^{i+1}0^{k+1}.
$$
Consider coefficients of the polynomial $V_k(y+1)$
$$
(\Delta^i0^k+\Delta^{i+1}0^k)(i+1)=
\left[\sum_{\alpha=0}^i(-1)^{i-\alpha}{i\choose \alpha}\alpha^k+
\sum_{\alpha=0}^{i+1}(-1)^{i+1-\alpha}{i+1\choose
\alpha}\alpha^k\right](i+1)=
$$
$$
=\left[\sum_{\alpha=0}^i(-1)^{i-\alpha}{i\choose \alpha}\alpha^k+
\sum_{\alpha=0}^i(-1)^{i+1-\alpha}{i+1\choose
\alpha}\alpha^k+(i+1)^k\right](i+1)=
$$
$$
=\left[\sum_{\alpha=0}^i(-1)^{i-\alpha}{i\choose \alpha}\alpha^k+
\sum_{\alpha=0}^i(-1)^{i+1-\alpha}\left({i\choose \alpha}+
{i\choose \alpha-1}\right)\alpha^k+(i+1)^k\right](i+1)=
$$
$$
=\left[\sum_{\alpha=0}^i(-1)^{i+1-\alpha}{i\choose
\alpha-1}\alpha^k+ (i+1)^k\right](i+1)=
$$
$$=
\sum_{\alpha=0}^{i+1}(-1)^{i+1-\alpha}{i+1\choose
\alpha}\alpha^{k+1} =\Delta^{i+1}0^{k+1}.
$$
Hence we get $V_k(y+1)=P_k(y+1)$, i.e.
$$
P_k(x)=E_k(x).
$$
Lemma 3 is proved.

\textbf{Proof of lemma 4.} Using (14) we have
$$
[D_m[\beta],[\beta]^{2m}]=(2m-1)!h^{-2m}\Bigg(2\sum\limits_{\beta=1}^{\infty}
\sum\limits_{k=1}^{m-1}\frac{(1-q_k)^{2m+1}}{q_kE_{2m-1}(q_k)}q_k^{\beta}[\beta]^{2m}+
$$
$$
+2[1]^{2m}\Bigg)=2\cdot(2m-1)!\left(\sum\limits_{\beta=1}^{\infty}\frac{(1-q_k)^{2m+1}}
{q_kE_{2m-1}(q_k)}q_k^{\beta}\cdot \beta^{2m}+1\right).
$$
Hence, by virtue of well-known formula
$\sum\limits_{\gamma=0}^{\infty}q^\gamma\gamma^k=\frac{1}{1-q}\sum_{i=0}^k
\left(\frac{q}{1-q}\right)^i\Delta^i0^k$, we obtain
$$
[D_m[\beta],[\beta]^{2m}]=2(2m-1)!\left(\sum\limits_{k=1}^{m-1}\frac{(1-q_k)^{2m+1}}
{q_kE_{2m-1}(q_k)}\cdot
\frac{1}{1-q_k}\sum\limits_{i=0}^{2m}\left(\frac{q_k}{1-q_k}\right)^i\Delta^i0^{2m}+1\right).
$$
According to lemma 3
$$
[D_m[\beta],[\beta]^{2m}]=2(2m-1)!\left(\sum\limits_{k=1}^{m-1}\frac{(1-q_k)^{2m+1}}
{q_kE_{2m-1}(q_k)}\cdot
\frac{E_{2m-1}(1/q_k)}{(1-q_k)(1/q_k-1)^{2m}}+1\right).
$$
Now, using lemma 2, finally we get
$$
[D_m[\beta],[\beta]^{2m}]=2(2m-1)!\left(\sum\limits_{k=1}^{m-1}\frac{(1-q_k)^{2m+1}}
{q_kE_{2m-1}(q_k)}\cdot
\frac{E_{2m-1}(q_k)}{(1-q_k)q_k^{2m-1}(1/q_k-1)^{2m}}+1\right)=
$$
$$
=2\cdot (2m-1)!\cdot m=(2m)!.
$$
Lemma 4 is proved.\\[0.2cm]

\textbf{Proof of theorem 3.}

Suppose $\beta=\overline{1,N-1}$. Then from (17), using (24) and
definition of convolution of discrete functions, we have
$$
C[\beta]=hD_m[\beta]*u[\beta]=h\sum\limits_{\gamma=-\infty}^{\infty}
D_m[\beta-\gamma]u[\gamma]=h\Bigg(\sum\limits_{\gamma=-\infty}^{-1}
D_m[\beta-\gamma]\left(-Q_{2m-1}[\gamma]+R_{m-1}^-[\gamma]\right)+
$$
$$
+\sum\limits_{\gamma=0}^ND_m[\beta-\gamma]f_m[\gamma]+
\sum\limits_{\gamma=N+1}^{\infty}D_m[\beta-\gamma]\left(Q_{2m-1}[\gamma]+R_{m-1}^+[\gamma]\right)\Bigg)=
$$
$$
=h\Bigg(\sum\limits_{\gamma=-\infty}^{\infty}D_m[\beta-\gamma]f_m[\gamma]+
\sum\limits_{\gamma=1}^{\infty}D_m[\beta+\gamma](-Q_{2m-1}[-\gamma]+
R_{m-1}^-[-\gamma]-f_m[-\gamma])+
$$
$$
+\sum\limits_{\gamma=1}^{\infty}D_m[N+\gamma-\beta](Q_{2m-1}[N+\gamma]+
R_{m-1}^+[N+\gamma]-f_m[N+\gamma])\Bigg)
$$
$$
=h\Bigg(D_m[\beta]*f_m[\beta]+
\sum\limits_{\gamma=1}^{\infty}D_m[\beta+\gamma](-Q_{2m-1}[-\gamma]+
R_{m-1}^-[-\gamma]-f_m[-\gamma])+
$$
$$
+\sum\limits_{\gamma=1}^{\infty}D_m[N+\gamma-\beta](Q_{2m-1}[N+\gamma]+
R_{m-1}^+[N+\gamma]-f_m[N+\gamma])\Bigg).
$$
Hence, using (14), (15) and (33), we get
$$
C[\beta]=h\Bigg(1+\sum\limits_{k=1}^{m-1}q_k^{\beta}\frac{(2m-1)!}{h^{2m}}\frac{(1-q_k)^{2m+1}}
{q_kE_{2m-1}(q_k)}\sum\limits_{\gamma=1}^{\infty}q_k^{\gamma}
\left(-Q_{2m-1}[-\gamma]+R_{m-1}^-[-\gamma]-f_m[-\gamma]\right)+
$$
$$
+\sum\limits_{k=1}^{m-1}q_k^{N-\beta}\frac{(2m-1)!}{h^{2m}}\frac{(1-q_k)^{2m+1}}
{q_kE_{2m-1}(q_k)}\sum\limits_{\gamma=1}^{\infty}q_k^{\gamma}
\left(Q_{2m-1}[N+\gamma]+R_{m-1}^+[N+\gamma]-f_m[N+\gamma]\right)\Bigg).\eqno
(35)
$$
We denote
$$
\begin{array}{l}
d_k=\frac{(2m-1)!}{h^{2m}}\frac{(1-q_k)^{2m+1}}
{q_kE_{2m-1}(q_k)}\sum\limits_{\gamma=1}^{\infty}q_k^{\gamma}
\left(-Q_{2m-1}[-\gamma]+R_{m-1}^-[-\gamma]-f_m[-\gamma]\right),\\
p_k=\frac{(2m-1)!}{h^{2m}}\frac{(1-q_k)^{2m+1}}
{q_kE_{2m-1}(q_k)}\sum\limits_{\gamma=1}^{\infty}q_k^{\gamma}
\left(Q_{2m-1}[N+\gamma]+R_{m-1}^+[N+\gamma]-f_m[N+\gamma]\right),
\end{array}\eqno (36)
$$
where $k=\overline{1,m-1}$. Then when $\beta=\overline{1,N-1}$ for
the optimal coefficients
$$
C[\beta]=h\left(1+\sum\limits_{k=1}^{m-1}\left(d_kq_k^{\beta}+p_kq_k^{N-\beta}\right)\right).
\eqno (37)
$$
Now from (8) when $\alpha=0,1$ keeping in mind (37) for $C[0]$ and
$C[N]$ we obtain the expressions
$$
C[0]=h\left(\frac{1}{2}+\sum_{k=1}^{m-1}\left(\frac{d_kq_k-p_kq_k^N}{q_k-1}-
\frac{h(q_k^{N+1}-q_k)(d_k-p_k)}{(q_k-1)^2}\right)\right), \eqno
(38)
$$
$$
C[N]=h\left(\frac{1}{2}+\sum_{k=1}^{m-1}\left(\frac{p_kq_k-d_kq_k^N}{q_k-1}+
\frac{h(q_k^{N+1}-q_k)(d_k-p_k)}{(q_k-1)^2}\right)\right), \eqno
(39)
$$
where $d_k$ and $p_k$ are determined from (36).

Now in order to prove theorem 3 sufficiently to show the
equalities $d_k=p_k$, $k=\overline{1,m-1}$.

Now we will prove, that $d_k=p_k$, $k=\overline{1,m-1}$.

Consider the convolution in equality (6) and we rewrite it in the
form
$$
g[\beta]=G_{m,1}[\beta]*C[\beta]=\sum_{\gamma=0}^NC[\gamma]\frac{|h\beta-h\gamma|^{2m-1}}{2(2m-1)!}=
$$
$$=
\sum_{\gamma=0}^{\beta}C[\gamma]\frac{(h\beta-h\gamma)^{2m-1}}{(2m-1)!}-
\sum_{\gamma=0}^{N}C[\gamma]\frac{(h\beta-h\gamma)^{2m-1}}{2(2m-1)!}=
 T_1-T_2,\eqno (40)
$$
where
$$
T_1=\sum_{\gamma=0}^{\beta}C[\gamma]\frac{(h\beta-h\gamma)^{2m-1}}{(2m-1)!},
\ \
T_2=\sum_{\gamma=0}^{N}C[\gamma]\frac{(h\beta-h\gamma)^{2m-1}}{2(2m-1)!}.
\eqno (41)
$$

First we consider $T_1$. For $T_1$ using (34) we have
$$
T_1=C[0]\frac{(h\beta)^{2m-1}}{(2m-1)!}+\sum\limits_{\gamma=1}^{\beta}h
\left(1+\sum\limits_{k=1}^{m-1}\left(d_kq_k^{\gamma}+
p_kq_k^{N-\gamma}\right)\right)\frac{(h\beta-h\gamma)^{2m-1}}{(2m-1)!}=
$$
$$
=C[0]\frac{(h\beta)^{2m-1}}{(2m-1)!}+\frac{h^{2m}}{(2m-1)!}\Bigg(
\sum_{\gamma=0}^{\beta-1}\gamma^{2m-1}+$$
$$+\sum\limits_{k=1}^{m-1}\left(
d_kq_k^{\beta}\sum_{\gamma=0}^{\beta-1}q_k^{-\gamma}\gamma^{2m-1}+
p_kq_k^{N-\beta}\sum_{\gamma=0}^{\beta-1}q_k^{\gamma}\gamma^{2m-1}
\right)\Bigg).
$$
Using by formulas (25), (26) the expression for $T_1$ we reduce to
the form
$$
T_1=C[0]\frac{(h\beta)^{2m-1}}{(2m-1)!}+\frac{h^{2m}}{(2m-1)!}\Bigg[
\sum\limits_{j=1}^{2m}\frac{(2m-1)!B_{2m-j}}{J!(2m-j)!}\beta^j+
$$
$$
+\sum\limits_{k=1}^{m-1}\Bigg(d_kq_k^{\beta}\left[\frac{q_k}{q_k-1}
\sum\limits_{i=0}^{2m-1}\frac{\Delta^i0^{2m-1}}{(q_k-1)^i}-\frac{q_k^{1-\beta}}{q_k-1}
\sum\limits_{i=0}^{2m-1}\frac{\Delta^i\beta^{2m-1}}{(q_k-1)^i}
\right]+
$$
$$
+p_kq_k^{N-\beta}\left[\frac{1}{1-q_k}\sum\limits_{i=0}^{2m-1}
\left(\frac{q_k}{1-q_k}\right)^i\Delta^i0^{2m-1}-\frac{q_k^{\beta}}{1-q_k}
\sum\limits_{i=0}^{2m-1}\left(\frac{q_k}{1-q_k}\right)^i\Delta^i\beta^{2m-1}\right]
\Bigg)\Bigg]\eqno (42)
$$
Keeping in mind that $q_k$ is the root of the Euler-Frobenius
polynomial $E_{2m-2}(x)$ of degree $2m-2$, using lemma 3 and
equality (38), after simplifications, from (42) we get
$$
T_1=\frac{h^{2m}}{(2m-1)!}\Bigg[\frac{\beta^{2m}}{2m}+\beta^{2m-1}
\frac{h(q_k^{N+1}-q_k)(p_k-d_k)}{(q_k-1)^2}+
$$
$$
+\sum\limits_{j=1}^{2m-2}\frac{(2m-1)!B_{2m-j}}{j!(2m-j)!}\beta^j+
\sum\limits_{k=1}^{m-1}\sum\limits_{i=1}^{2m-1}\frac{-d_kq_k+(-1)^iq_k^{N+i}
p_k}{(q_k-1)^{i+1}}\Delta^i\beta^{2m-1}\Bigg].
$$
Hence using the formula
$$
\Delta^i\beta^{2m-1}=\sum\limits_{j=0}^{2m-1}{2m-1 \choose
j}\Delta^i0^j\beta^{2m-1-j}
$$
and grouping in powers of $\beta$ we have
$$
T_1=\frac{(h\beta)^{2m}}{(2m)!}+\frac{h^{2m+1}}{(2m-1)!}\beta^{2m-1}
\frac{(q_k^{N+1}-q_k)(p_k-d_k)}{(q_k-1)^2}+
$$
$$
+\frac{h^{2m}}{(2m-1)!}\sum\limits_{j=1}^{2m-2}\frac{(2m-1)!\beta^{2m-1-j}}
{j!(2m-1-j)!}\left[\frac{B_{j+1}}{j+1}+\sum\limits_{k=1}^{m-1}\sum\limits_{i=1}^j
\frac{-d_kq_k+(-1)^iq_k^{N+i}p_k}{(q_k-1)^{i+1}}\Delta^i0^j\right]+
$$
$$
+\frac{h^{2m}}{(2m-1)!}\sum\limits_{k=1}^{m-1}\sum\limits_{i=1}^{2m-1}
\frac{-d_kq_k+(-1)^iq_k^{N+i}p_k}{(q_k-1)^{i+1}}\Delta^i0^{2m-1}.\eqno
(43)
$$

Now consider $T_2$. Using binomial formula and equalities (5) the
expression for $T_2$ we reduce to the form
$$
T_2=\sum\limits_{\gamma=0}^NC[\gamma]\frac{(h\beta-h\gamma)^{2m-1}}{2(2m-1)!}=
$$
$$
=\frac{1}{2}\sum\limits_{j=0}^{m-1}\frac{(h\beta)^{2m-1-j}(-1)^j}{(j+1)!
(2m-1-j)!}+\frac{1}{2}\sum\limits_{j=1}^m\frac{(h\beta)^{m-j}(-1)^{m+j-1}}
{(m+j-1)!(m-j)!}\sum\limits_{\gamma=0}^NC[\gamma][\gamma]^{m+j-1}.\eqno
(44)
$$
And for the right hand side $f_m[\beta]$ of equation (4) (see
(10)) using binomial formula we get
$$
f_m[\beta]=\frac{(h\beta)^{2m}}{(2m)!}+\frac{1}{2}\sum\limits_{j=0}^{m-1}
\frac{(h\beta)^{2m-1-j}(-1)^{j+1}}{(j+1)!(2m-1-j)!}+\sum\limits_{j=1}^{m}
\frac{(h\beta)^{m-j}(-1)^{m+j}}{(m+j)!(m-j)!}.\eqno (45)
$$
Then substituting (43), (44) into (40) taking into account (45)
after some simplifications for the difference
$f_m[\beta]-g[\beta]$ we get
$$
f_m[\beta]-g[\beta]=\frac{h^{2m+1}}{(2m-1)!}\beta^{2m-1}\sum\limits_{k=1}^{m-1}
\frac{(q_k^{N+1}-q_k)(d_k-p_k)}{(q_k-1)^2}-
$$
$$
-\frac{h^{2m}}{(2m-1)!}\sum\limits_{j=1}^{2m-2}\frac{(2m-1)!\beta^{2m-1-j}}
{j!(2m-1-j)!}\left[\frac{B_{j+1}}{j+1}+\sum\limits_{k=1}^{m-1}\sum\limits_{i=1}^j
\frac{-d_kq_k+(-1)^iq_k^{N+i}p_k}{(q_k-1)^{i+1}}\Delta^i0^j\right]+
$$
$$
+\frac{1}{2}\sum\limits_{j=1}^m\frac{(h\beta)^{m-j}(-1)^{m+j}}{(m+j-1)!(m-j)!}\left[
\frac{1}{m+j}-\sum\limits_{\gamma=0}^NC[\gamma][\gamma]^{m+j-1}\right]-
$$
$$
-\frac{h^{2m}}{(2m-1)!}\sum\limits_{k=1}^{m-1}\sum\limits_{i=1}^{2m-1}
\frac{-d_kq_k+(-1)^iq_k^{N+i}p_k}{(q_k-1)^{i+1}}\Delta^i0^{2m-1}.
$$
Hence note that $f_m[\beta]-g[\beta]$ is polynomial of degree
$(2m-1)$ of $[\beta]=h\beta$, i.e.
$$
f_m[\beta]-g[\beta]=\sum\limits_{j=0}^{2m-1}a_j[\beta]^j.\eqno
(46)
$$
Here
$$
a_j= \left\{
\begin{array}{ll}
\frac{h^{2}}{(2m-1)!}\sum\limits_{k=1}^{m-1}
\frac{(q_k^{N+1}-q_k)(d_k-p_k)}{(q_k-1)^2}& \mbox{ for } j=2m-1,\\
b_j & \mbox{ for } m\le j\leq 2m-2,\\
b_j+\frac{(-1)^j}{2j!(2m-1-j)!}\left(\frac{1}{2m-j}-\sum\limits_{\gamma=0}^{N}
C[\gamma][\gamma]^{2m-j-1}\right)& \mbox{ for } 1\leq j\leq m-1,\\
-\frac{h^{2m}}{(2m-1)!}\sum\limits_{k=1}^{m-1}\sum\limits_{i=1}^{2m-1}
\frac{-d_kq_k+(-1)^iq_k^{N+i}p_k}{(q_k-1)^{i+1}}\Delta^i0^{2m-1}+&\\
+\frac{1}{2(2m-1)!}\left(\frac{1}{2m}-\sum\limits_{\gamma=0}^{N}
C[\gamma][\gamma]^{2m-1}\right)& \mbox{ for } j=0,\\
\end{array}
\right.
$$
where
$$
b_j=-\frac{h^{2m-j}}{j!(2m-j-1)!}\left(\frac{B_{2m-j}}{2m-j}+
\sum\limits_{i=1}^{2m-j-1}\sum\limits_{k=1}^{m-1}\frac{-d_k
q_k+(-1)^iq_k^{N+i}p_k}{(q_k-1)^{i+1}}\Delta^i0^{2m-j-1}\right).
$$
On the other hand from (6) for difference $f_m[\beta]-g[\beta]$
the following is true
$$
f_m[\beta]-g[\beta]=P_{m-1}[\beta].\eqno (47)
$$
If
$$
\sum\limits_{k=1}^{m-1}\sum_{i=1}^j\frac{d_kq_k+(-1)^{i+1}q_k^{N+i}p_k}
{(q_k-1)^{i+1}}\Delta^i0^j=\frac{B_{j+1}}{j+1},\ \
j=\overline{1,m-1},\eqno (48)
$$
$$
\sum\limits_{k=1}^{m-1}
\frac{(q_k^{N+1}-q_k)(d_k-p_k)}{(q_k-1)^2}=0,\eqno (49)
$$
then equality (47) is take placed.  From (46) and (47) we find
unknown polynomial $P_{m-1}[\beta]$ in the system (6)-(8)
$$
P_{m-1}[\beta]=\sum\limits_{j=1}^{m-1}a_j[\beta]^j.\eqno (50)
$$
Thus, for unknowns $d_k$ and $p_k$ ($k=\overline{1,m-1}$) we have
got system of linear equations (48) -(49).

Later, from (8) for $\alpha=\overline{2,m-1}$, using equalities
(37), (39), (25), (26) and (49) after some simplifications we
obtain
$$
h^{\alpha+1}\sum\limits_{k=1}^{m-1}\sum\limits_{i=0}^{\alpha}
\left(\frac{d_kq_k^i+(-1)^{i+1}q_k^{N+1}p_k}{(1-q_k)^{i+1}}
\Delta^i0^{\alpha}-\frac{d_kq_k^{N+i}+
(-1)^{i+1}q_kp_k}{(1-q_k)^{i+1}}\Delta^i0^{\alpha}\right)+
$$
$$
+\sum\limits_{j=2}^{\alpha}\frac{\alpha!h^j}{(j-1)!(\alpha+1-j)!}
\left[\frac{B_j}{j}-\sum\limits_{k=1}^{m-1}\sum\limits_{i=0}^{j-1}
\frac{d_kq_k^{N+i}+(-1)^{i+1}q_kp_k}{(1-q_k)^{i+1}}\Delta^i0^{j-1}\right]=0.
\eqno (51)
$$
The left hand side of the equation (51) is the polynomial of
degree $(\alpha+1)$ with respect to $h$. Then from (51) obviously
that all coefficients of this polynomial are zero, i.e.
$$
\sum\limits_{k=1}^{m-1}\sum\limits_{i=0}^{j}\frac{d_kq_k^{N+i}+(-1)^{i+1}q_kp_k}
{(1-q_k)^{i+1}}\Delta^i0^j=\frac{B_{j+1}}{j+1},\ \
j=\overline{1,\alpha-1},\eqno (52)
$$
$$
\sum_{k=1}^{m-1}\sum_{i=0}^{\alpha}\frac{d_kq_k^i+(-1)^{i+1}q_k^{N+1}p_k}
{(1-q_k)^{i+1}}\Delta^i0^{\alpha}=\sum_{k=1}^{m-1}\sum_{i=0}^{\alpha}
\frac{d_kq_k^{N+i}+(-1)^{i+1}q_kp_k}{(1-q_k)^{i+1}}\Delta^i0^{\alpha},\
\ \alpha=\overline{2,m-1}. \eqno (53)
$$
Thus for unknowns $d_k$ and $p_k$ we have got the system (52),
(53).

Using lemmas 2, 3 for the left hand side of the system (53) we
obtain
$$
\sum_{k=1}^{m-1}\sum_{i=0}^{\alpha}\frac{d_kq_k^i+(-1)^{i+1}q_k^{N+1}p_k}
{(1-q_k)^{i+1}}\Delta^i0^{\alpha}=(-1)^{\alpha+1}
\sum_{k=1}^{m-1}\sum_{i=0}^{\alpha}\frac{d_kq_k+(-1)^{i+1}q_k^{N+i}p_k}
{(q_k-1)^{i+1}}\Delta^i0^{\alpha},
$$
where $\alpha=2,3,...,m-1$.\\
Hence taking into account (48) we get
$$
\sum_{k=1}^{m-1}\sum_{i=0}^{\alpha}\frac{d_kq_k^i+(-1)^{i+1}q_k^{N+1}p_k}
{(1-q_k)^{i+1}}\Delta^i0^{\alpha}=(-1)^{\alpha+1}\frac{B_{\alpha+1}}{\alpha+1},
\ \ \alpha=2,3,...,m-1.\eqno (54)
$$
Since for $\alpha=2\ell,\ \ell=1,2,...$ Bernoulli numbers
$B_{\alpha+1}=B_{2\ell+1}$ are zero, then combining (54) and (53)
we have
$$
\sum_{k=1}^{m-1}\sum_{i=0}^{\alpha}
\frac{d_kq_k^{N+i}+(-1)^{i+1}q_kp_k}{(1-q_k)^{i+1}}\Delta^i0^{\alpha}
=\frac{B_{\alpha+1}}{\alpha+1}, \ \ \alpha=2,3,...,m-1.\eqno (55)
$$
Comparing systems (52) and (55) it is easy to see, that from (52)
for $j=1$ we get only one equation with respect to unknowns $d_k$
and $p_k$. And for $j=2,...,\alpha-1$ obtained equations from
system (52) are part of the system (55).

Thus, from (48), (52), (55) we get the following system
$$
\sum\limits_{k=1}^{m-1}\sum_{i=1}^j\frac{d_kq_k+(-1)^{i+1}q_k^{N+i}p_k}
{(q_k-1)^{i+1}}\Delta^i0^j=\frac{B_{j+1}}{j+1},\ \
j=1,2,...,m-1,\eqno (56)
$$
$$
\sum_{k=1}^{m-1}\sum_{i=0}^{j}
\frac{d_kq_k^{N+i}(-1)^{i+1}+q_kp_k}{(q_k-1)^{i+1}}\Delta^i0^{j}
=\frac{B_{j+1}}{j+1}, \ \ j=1,2,...,m-1.\eqno (57)
$$
Subtracting equality (57) from (56) gives
$$
\sum_{k=1}^{m-1}\sum_{i=0}^{j}\frac{(q_k-q_k^{N+i}(-1)^{i+1})(d_k-p_k)}
{(q_k-1)^{i+1}}\Delta^i0^j=0,\ \ j=1,2,...,m-1.\eqno (58)
$$
From system (58) for $j=1$ we get equation (49). Therefore for
finding unknowns $d_k$ and $p_k$ it is sufficiently to solve the
system (56), (57).

Further we investigate the system (58). After some transformations
homogenous system (58) we reduce to the form
$$
\sum_{k=1}^{m-1}\frac{(q_k-q_k^{N+j}(-1)^{j+1})(d_k-p_k)}
{(q_k-1)^{j+1}}=0,\ \ j=1,2,...,m-1,\eqno (59)
$$
where $d_k-p_k$, $k=1,2,...,m-1$ are unknowns. Obviously, that the
main matrix of the system (59) is the Vandermonde type matrix,
i.e. system (59) has a unique solution, which identically zero.
This means
$$
d_k=p_k,\ \ k=1,2,...,m-1.\eqno (60)
$$
Then keeping in mind (60) from (37), (38), (39) and (56) we
respectively get the following system
$$
C[\beta]=h\left(1+\sum_{k=1}^{m-1}d_k\left(q_k^{\beta}+q_k^{N-\beta}\right)\right),\
\ \beta=\overline{1,N-1},
$$
$$
C[0]=h\left(\frac{1}{2}-\sum_{k=1}^{m-1}d_k\frac{q_k-q_k^N}{1-q_k}\right),
$$
$$
C[N]=h\left(\frac{1}{2}-\sum_{k=1}^{m-1}d_k\frac{q_k-q_k^N}{1-q_k}\right),
$$
$$
\sum\limits_{k=1}^{m-1}d_k\sum_{i=1}^j\frac{q_k+(-1)^{i+1}q_k^{N+i}}
{(q_k-1)^{i+1}}\Delta^i0^j=\frac{B_{j+1}}{j+1},\ \ j=1,2,...,m-1.
$$

Theorem 3 is proved.

\end{document}